


\magnification\magstephalf

\hoffset=.5truein \advance\hoffset by -5.5truemm
\hsize=5.7truein
\voffset=.4375truein
\vsize=8.25truein

\parskip=6pt plus2pt minus1pt
\smallskipamount=6pt plus2pt minus2pt
\medskipamount=12pt plus4pt minus4pt
\bigskipamount=18pt plus6pt minus6pt
\abovedisplayskip=12pt plus4pt minus4pt
\belowdisplayskip=12pt plus4pt minus4pt

\normalbaselineskip=14pt plus1pt minus1pt
\normallineskiplimit=1pt
\normallineskip=1pt
\normalbaselines


\font\titlefont    = cmbx12 at 16pt
\font\twelverm     = cmr12
\font\ninerm       = cmr9
\font\sixrm        = cmr6
\font\nineit       = cmti9
\font\twelvebf     = cmbx12
\font\ninebf       = cmbx9
\font\eightbf      = cmbx8
\font\sixbf        = cmbx6

\font\twelvebfi    = cmmib10 at 12pt
\font\eightbfi     = cmmib8
\font\ninei        = cmmi9
\font\sixi         = cmmi6
\font\twelvebfsy   = cmbsy10 at 12pt
\font\eightbfsy    = cmbsy8
\font\ninesy       = cmsy9      
\font\sixsy        = cmsy6

\font\tenmsa       = msam10
\font\ninemsa      = msam9
\font\sevenmsa     = msam7
\font\sixmsa       = msam6
\font\tenmsb       = msbm10
\font\ninemsb      = msbm9
\font\sevenmsb     = msbm7
\font\sixmsb       = msbm6
\font\teneuf       = eufm10   
\font\nineeuf      = eufm9
\font\seveneuf     = eufm7
\font\sixeuf       = eufm6
\font\teneus       = eusm10   
\font\twelveeusb   = eusb10 at 12pt 
\font\nineeus      = eusm9
\font\seveneus     = eusm7
\font\sixeus       = eusm6


\newfam\msafam
\textfont\msafam           =\tenmsa 
\scriptfont\msafam         =\sevenmsa

\newfam\msbfam
\def\bbb{\fam\msbfam}
\textfont\msbfam           =\tenmsb 
\scriptfont\msbfam         =\sevenmsb

\newfam\euffam
\def\euf{\fam\euffam}
\textfont\euffam           =\teneuf 
\scriptfont\euffam         =\seveneuf 

\newfam\eusfam
\def\eus{\fam\eusfam}
\textfont\eusfam           =\teneus 
\scriptfont\eusfam         =\seveneus

\def\ninepoint{
  \textfont0                 = \ninerm
  \scriptfont0               = \sixrm
  \textfont1                 = \ninei
  \scriptfont1               = \sixi
  \textfont2                 = \ninesy 
  \scriptfont2               = \sixsy
  \textfont\bffam            = \ninebf   
  \scriptfont\bffam          = \sixbf   
  \textfont\msafam           = \ninemsa
  \scriptfont\msafam         = \sixmsa
  \textfont\msbfam           = \ninemsb
  \scriptfont\msbfam         = \sixmsb
  \textfont\euffam           = \nineeuf
  \scriptfont\euffam         = \sixeuf
  \textfont\eusfam           = \nineeus
  \scriptfont\eusfam         = \sixeus
  \def\rm{\ninerm\fam0}%
  \def\it{\nineit\fam1}%
  \def\bf{\ninebf\fam\bffam}%
  \rm
  \normalbaselineskip=12pt%
  \normallineskiplimit=1pt%
  \normallineskip=1pt%
  \normalbaselines}

\def\boldtwelvepoint{
  \textfont0                 = \twelvebf
  \scriptfont0                 = \eightbf
  \textfont1                 = \twelvebfi
  \scriptfont1                = \eightbfi
  \textfont2                 = \twelvebfsy
  \scriptfont2                 = \eightbfsy
  \textfont\eusfam           = \twelveeusb
  \def\rm{\twelvebf\fam0}%
  \normalbaselineskip=18pt%
  \normallineskiplimit=0pt%
  \normallineskip=1.2pt%
  \normalbaselines\rm}


\headline={\ifnum\pageno>1\ifodd\pageno\hfil{\ninerm\folio}%
  \else{\ninerm\folio}\hfil\fi\else\hfil\fi}
\footline={\hfil}
\def\makeheadline{\vbox to0pt{\vskip-34pt%
   \line{\vbox to8.5pt{}\the\headline}\vss}\nointerlineskip}

\def\titlelines#1\par{\leavevmode\bigskip\vskip\parskip\nointerlineskip
  \centerline{\vbox{\baselineskip=24pt\titlefont
  \halign{\hfil##\hfil\cr#1\crcr}}}\par}
\def\authorlines#1\par{\medskip\vskip\parskip\nointerlineskip
  \centerline{\vbox{\twelverm\baselineskip=18pt
  \halign{\hfil##\hfil\cr#1\crcr}}}\par}

\newif\ifheadingfeed
\def\heading#1\par{%
  \ifheadingfeed\vskip0pt plus.3\vsize\penalty-250\vskip0pt plus-.3\vsize
    \else\headingfeedtrue\fi
  \bigskip\vskip\parskip
  \centerline{\vbox{\boldtwelvepoint
  \halign{\hfil##\hfil\cr#1\crcr}}}%
  \smallskip\noindent}
\def\subheading#1\par{\removelastskip\smallskip
  \penalty-75\noindent{\boldtenpoint#1}\enspace}
\def\subsubheading#1\par{\penalty-50{\it#1}\enspace}
\long\def\proclaim#1\par#2\par{\ifdim\lastskip<\smallskipamount
  \removelastskip\smallskip\fi\penalty-100{\bf#1}\enspace
  {\parskip=3pt\it#2}\smallskip\penalty100}
\def\demo#1\par{\ifdim\lastskip<\smallskipamount\removelastskip\smallskip
  \penalty-50\fi{\it#1}\enspace}
\def\enddemo{\penalty200\hbox{}\nobreak\hfill\space\hbox{$\square$}\smallskip
  \penalty-75}
\def\item#1\par#2\par{{\advance\leftskip by 3\parindent
  \noindent\kern-2\parindent\rlap{#1}\kern2\parindent#2\par}}
\def\eqalign#1{\null\,\vcenter{\openup\jot
  \ialign{&\strut\hfil$\displaystyle##$%
          &$\displaystyle{}##$\hfil\cr#1\crcr}}\,}
\catcode`\@=11
\def\eqalignno#1{\displ@y\tabskip=\centering
  \halign to\displaywidth{\hfil$\@lign\displaystyle{##}$\tabskip=0pt
    &$\@lign\displaystyle{{}##}$\hfil\tabskip=\centering
    &\kern-\displaywidth\rlap{\hskip\parindent$\@lign##$}%
    \tabskip=\displaywidth\crcr#1\crcr}}
\def\displaylinesno#1{\displ@y\tabskip=\centering
  \halign to\displaywidth{\hfil$\@lign\displaystyle{##}$\hfil
  &\kern-\displaywidth\rlap{\hskip\parindent$\@lign##$}%
  \tabskip=\displaywidth\crcr#1\crcr}}
\catcode`\@=12
\def\references{\heading\tenbf References\par\smallskip
  \parskip=2pt plus.7pt minus.7pt\tolerance=1000\overfullrule=0pt
  \frenchspacing\ninepoint}
\def\ref[#1]#2\par {\noindent\kern-1.8cm\hbox to0.8cm{[#1]\hfil}#2%
  \leftskip= 1.8cm\par\leftskip= 0cm}
\def\endlines#1{\medskip\vskip\parskip\nointerlineskip
  \centerline{\vbox{\hrule width 1in}}
  \medskip\nointerlineskip\vskip\parskip
  \centerline{\ninepoint\vbox{\halign{\hfil##\hfil\cr#1\crcr}}}}
\def\date{\ifcase\month\or January\or February\or March\or April\or May%
  \or June\or July\or August\or September\or October\or November\or December%
  \fi\space\number\day,\space\number\year}

\newif\ifshowlabel\showlabelfalse
\def\eqno#1{\leqno{%
  \ifshowlabel\llap{\hbox to.5in{\eightrm\string#1\hss}}\fi
  \hskip\parindent#1}}
\catcode`\@=11
\def\eqalignno#1{\displ@y\tabskip=\centering
  \halign to\displaywidth{\hfil$\@lign\displaystyle{##}$\tabskip=0pt
    &$\@lign\displaystyle{{}##}$\hfil\tabskip=\centering
    &\kern-\displaywidth\rlap{\hskip\parindent$\@lign##$}%
    \tabskip=\displaywidth\crcr#1\crcr}}
\catcode`\@=12

\hyphenation{mani-fold mani-folds eigen-space eigen-vector}


\def\Ad{\mathop{\rm Ad}}
\def\ad{\mathop{\rm ad}}
\def\cont{{\rm cont}}
\def\Hom{\mathop{\rm Hom}}
\def\ss{{\rm ss}}
\def\Tr{\mathop{\rm Tr}}

\def\U{{\rm U}}
\def\SO{{\rm SO}}
\def\SU{{\rm SU}}
\def\Spin{{\rm Spin}}

\def\IR{{\bbb R}}
\def\IZ{{\bbb Z}}

\def\eufg{{\euf g}}
\def\eufh{{\euf h}}
\def\eufk{{\euf k}}

\def\eufz{{\euf z}}
\def\eufhom{\mathop{\euf hom}}
\def\eufX{{\euf X}}

\def\eusA{{\eus A}}
\def\eusG{{\eus G}}
\def\eusGhat{\widehat{\eus G}}
\def\eusO{{\eus O}}

\def\epsilon{\varepsilon}
\def\phi{\varphi}

\mathchardef\square="0803
\mathchardef\boxtimes="2802

\def\nablaA{\nabla_{\!\!A}}
\def\mapright#1{\buildrel#1\over\longrightarrow}

\def\mi{d}
\def\YMH{{\eus YMH}}
\def\YM{{\eus YM}}


\titlelines Compactness Theorems for\cr Invariant Connections\cr

\authorlines Johan R\aa de

\medskip\smallskip

\centerline{\vbox{\hrule width 2in}}

\smallskip

%
%


\heading \S1. Condition~C for invariant Yang-Mills

The Palais-Smale Condition~C holds for the Yang-Mills functional
on principal bundles over compact manifolds of dimension $\le 3$.
This was established by S.~Sedlacek [17] and C.~Taubes [18] Proposition~4.5
using the compactness theorem of K.~Uhlenbeck [20];
see also [23].
It is well known that Condition~C fails for Yang-Mills over compact manifolds
of dimension $\ge 4$.
The example of $\SO(3)$-invariant $\SU(2)$-connections over $S^4$,
see [2], [14], and [16],
suggested that Condition~C holds
for Yang-Mills over compact manifolds of any dimension
when restricted to connections that are
invariant under a group action on the manifold 
with orbits of codimension $\le 3$.
Such a result, essentially Theorem~3 below, was established by T.~Parker [14].
In this paper we generalize his result.

Let $X$ be a smooth compact Riemannian manifold of any dimension,
let $G$ be a compact Lie group, and let $P$ be a smooth principal
$G$-bundle over $X$.
Let $\eusA$ denote the space of smooth connections on $P$
and $\eusG$ the group of smooth bundle automorphisms of $P$.
We denote the gauge equivalence class of $A\in\eusA$ by $[A]\in\eusA/\eusG$.
We fix a $G$-invariant positive definite 
inner product on the Lie algebra $\eufg$.
Then the Yang-Mills functional is given by 
$$\YM[A]={1\over2}\int_X|F_A|^2dx$$
where $F_A$ is the curvature of $A$.
The critical points of the Yang-Mills functional
are the solutions of the Yang-Mills equation $d_A^*F_A=0$.

Let $H$ be a compact Lie group.
We fix a smooth action $\rho:H\times X\to X$ by isometries.
We then consider the set of smooth lifts
of $\rho$ to actions $\sigma:H\times P\to P$ by bundle maps.
The group $\eusG$ of gauge transformations acts on this set;
$g.\sigma(h) = g\,\sigma(h)\,g^{-1}$.
We choose one representative 
$$\sigma_i:H\times P\to P$$
from each gauge equivalence class of lifts of $\rho$;
here $i$ ranges over some index set $I$.
The action $\sigma_i$ of $H$ on $P$ induces actions of $H$ on $\eusA$ and 
$\eusG$. Let $\eusA^H_i$ and $\eusG^H_i$ denote the fixed point sets 
of these actions.
The group $\eusG^H_i$ acts on the affine space $\eusA^H_i$ 
and we form the quotient $\eusA^H_i/\eusG^H_i$.
We denote the equivalence class of $A\in\eusA^H_i$ 
in $\eusA^H_i/\eusG^H_i$ by $[A,\sigma_i]$;
see \S3.
There are natural maps
$$\pi_i:\eusA^H_i/\eusG^H_i\to\eusA/\eusG$$
given by $\pi_i([A,\sigma_i])=[A]$.

{\it Example 1.}\enspace
Let $X=S^1$, $G=\U(1)$, and $H=\SO(2)$.
There is a unique principal $\U(1)$-bundle over $S^1$.
Let $\SO(2)$ act in the natural way on $S^1$.
This action has a unique lift, up to gauge equivalence.
Thus we omit the index $i$.
The space $\eusA^{\SO(2)}/\eusG^{\SO(2)}$ is a copy of $\IR$,
its image in $\eusA/\eusG$ is a copy of $S^1$,
and $\pi$ is the universal covering map.
The Yang-Mills functional of course vanishes.

{\it Example 2.}\enspace
Let $X=S^3$, $G=\SU(2)$, and $H=\Spin(4)$.
There is a unique principal $\SU(2)$-bundle over $S^3$.
Let $\Spin(4)$ act in the natural way on $S^3$.
This action has exactly two lifts, up to gauge equivalence,
the trivial lift $\sigma_0$ and a nontrivial lift $\sigma_1$.
Then $\eusA^{\Spin(4)}_0/\eusG^{\Spin(4)}_0$
consists of a single point while $\eusA^{\Spin(4)}_1/\eusG^{\Spin(4)}_1$
is a copy of $\IR$.
The image in $\eusA/\eusG$ of $\eusA^{\Spin(4)}_0/\eusG^{\Spin(4)}_0$ 
is the trivial connection,
while the image of $\eusA^{\Spin(4)}_1/\eusG^{\Spin(4)}_1$
is a curve in $\eusA/\eusG$ with a self intersection point 
at the trivial connection.
The Yang-Mills functional of course vanishes on 
$\eusA^{\Spin(4)}_0/\eusG^{\Spin(4)}_0$
while on $\eusA^{\Spin(4)}_1/\eusG^{\Spin(4)}_1$,
suitably identified with $\IR$, it is given by a constant times the polynomial
$(t^2-1)^2$.
The image in $\eusA/\eusG$ of the points $t=\pm1$
in $\eusA^{\Spin(4)}_1/\eusG^{\Spin(4)}_1$ 
is the trivial connection.
The image of the point $t=0$
is gauge equivalent to the Levi-Civita connection on $S^3$.

Next we discuss the Palais-Smale Condition~C for the Yang-Mills functional
for invariant connections.
The above examples show that the maps $\pi_i$ are not in general injective.
Hence there are two forms of Condition~C for invariant Yang-Mills.
In the following 
$\|\cdot\|_{L^{-1,2}_A}$ is the norm on $L^{-1,2}$
dual to the norm
$\|b\|_{L^{1,2}_A}^2
  = \|b\|_{L^2}^2
    + \|\nablaA b\|_{L^2}^2 
$
on $L^{1,2}$.
The Morrey space $L^{1,2}_3$ is discussed briefly in \S2 and at length 
in Appendix~A.
Note that on manifolds of dimension $\le 3$, $L^{1,2}_3=L^{1,2}$.

\proclaim Definition.

A sequence $A_k\in\eusA$ is said to be a Palais-Smale sequence
if there exists $M>0$ such that
$$\displaylines{
    \|F_{A_k}\|_{L^2(X,\Lambda^2T^*X\otimes\Ad P)}  < M
      \rlap{\qquad for all $k$} \cr
    \noalign{\line{and\hfill}}
    \big\|d^*_{A_k}F_{A_k}\big\|_{L^{-1,2}_{A_k}(X,T^*X\otimes\Ad P)}
       \to 0 \rlap{\qquad as $k\to \infty$.} \cr
          }
$$

\proclaim Definition.

Strong Condition~C is said to hold on $\eusA^H_i$
if there for any Palais-Smale sequence $A_k\in\eusA^H_i$
exists a subsequence, which we also denote $A_k$,
and a sequence $g_k\in\eusG^H_i$ such that
$g_k.A_k\in\eusA$ converges in $L^{1,2}_3$ to a Yang-Mills
connection $A_\infty\in\eusA^H_i$ as $k\to\infty$.

\proclaim Definition.

Weak Condition~C is said to hold on $\eusA^H_i$
if there for any Palais-Smale sequence $A_k\in\eusA^H_i$
exists a subsequence, which we also denote $A_k$,
and a sequence $g_k\in\eusG$ such that
$g_k.A_k\in\eusA$ converges in $L^{1,2}_3$ to a Yang-Mills
connection $A_\infty\in\eusA^H_i$\break as $k\to\infty$.

In strong Condition~C, $g_k\in\eusG^H_i$.
Thus strong Condition~C is essentially Condition~C for Yang-Mills
on $\eusA^H_i/\eusG^H_i$.
This space has the structure of an infinite dimensional orbifold,
similar to the structure of $\eusA/\eusG$.
If the action of $H$ on $X$ is free, then strong Condition~C 
for invariant Yang-Mills on $X$
is equivalent to Condition~C for the corresponding Yang-Mills-Higgs functional
on $X/H$.

In weak Condition~C, $g_k\in\eusG$.
Thus weak Condition~C is essentially Condition~C for Yaang-Mills
on the image 
$\pi_i(\eusA^H_i/\eusG^H_i)$
of $\eusA^H_i$ in $\eusA/\eusG$.
This can be a rather singular subset of $\eusA/\eusG$.
In the situation in the proof of Corollary~3,
the trivial connection is a self intersection point in 
$\pi(\eusA^{\SO(2)}/\eusG^{\SO(2)})$
of infinite multiplicity.

The following are our two main theorems.
They are proven in \S2.
Recall that $X$ can be of any dimension
and that $G$ and $H$ are only assumed to be compact.

\proclaim Theorem~1.

If all orbits of the action of $H$ on $X$ have codimension $\le 3$,
then weak Condition~C holds on $\eusA^H_i$ for all $i\in I$.

\proclaim Theorem~2.

Assume that all orbits of the action of $H$ on $X$ have codimension $\le 3$.
Then strong Condition~C holds on $\eusA^H_i$
if and only if every Yang-Mills connection
in $\eusA/\eusG$ has only finitely many preimages in $\eusA^H_i/\eusG^H_i$
under the map $\pi_i$.

The next two theorems give explicit 
sufficient conditions
for strong Condition~C.
They are special cases of Theorem~2, obtained in \S3 using Proposition~3.1.

\proclaim Theorem~3. \rm(T.~Parker [14])

If all orbits of the action of $H$ on $X$ have codimension $\le 3$,
$G$ is semisimple, and every Yang-Mills connection in $\eusA^H_i$
is irreducible,
then strong Condition~C holds on $\eusA^H_i$.

\proclaim Theorem~4.

If all orbits of the action of $H$ on $X$ have codimension $\le 3$
and $H$ is semisimple,
then strong Condition~C holds on $\eusA^H_i$ for all $i\in I$.

Theorem~3 was established by Parker on 4-manifolds
with the slightly stronger assumption that every connection in $\eusA^H_i$
is irreducible.
Under this assumption $\eusA^{H}_i/\eusG^{H}_i$ is a finite covering space
of $\pi_i(\eusA^{H}_i/\eusG^{H}_i)$.

\penalty-100
The following is a straightforward consequence of Theorem~1.

\proclaim Corollary~1.

If all orbits of the action of $H$ on $X$ have codimension $\le 3$,
then the Yang-Mills functional on $\eusA^H_i$
attains its infimum.

As an application we show how Theorem~4 gives a simple 
proof of the following result,
similar to the result of H.Y.~Wang [22].

\proclaim Corollary~2.

There exist irreducible non-(anti-)selfdual Yang-Mills\/ $\SU(2)$-con\-nections
over $S^2\times S^2$, with the standard metric,
of any even second Chern number not equal to $\pm2$ or $\pm4$.

This is proven in \S4 by applying ``the mountain pass lemma''
to paths between reducible Yang-Mills connections in
spaces of connections that are invariant 
under a $\Spin(3)$-action on $S^2\times S^2$.
H.-Y.~Wang [22] showed,
using a technically complicated gluing scheme,
that there exist irreducible $\SU(2)$ Yang-Mills 
connections over $S^2\times S^2$ of second Chern number $\pm2pq$, 
for any integers $p\ge q\ge 0$ such that $p\ne q(2k+1) + k(k+1)$ 
for all nonnegative integers $k$.
Thus Wang covered second Chern number $\pm4$, 
but not $\pm2$, $\pm10$, $\pm22,\,\ldots$
(The non-(anti-)selfdual $\SU(2)$ Yang-Mills connection
over $S^2\times S^2$ constructed in [21] Proposition~9.4 is reducible.)
For another application of Theorem~4, see U.~Gritsch [10].

Finally we give an application of Theorem~3.
This example shows that nontrivial phenomena
occur already in dimensions less than four.
Let $\YMH$ be the twisted Yang-Mills-Higgs functional
$$\YMH[A,\Phi]
  = {1\over2}\int_{S^2} \Bigl( \bigl|F_A-{*}\Phi\bigr|^2  
                               + \bigl|d_A\Phi\bigr|^2 
                        \Bigr)\,dx
$$
on $S^2$ with the standard metric.
Here $A$ is a connection on the trivial $\SU(2)$-bundle over $S^2$
and $\Phi$ is a section of the adjoint bundle.
The functional $\YMH$ over $S^2$ is the dimensional reduction
of the $\SU(2)$ Yang-Mills functional over $S^3$ under the
Hopf action of $\SO(2)$ on $S^3$.
The following result is proven in \S5.

\proclaim Corollary~3.

{
{\bf a)} The minimum of $\YMH$ is 0.
    There are countably infinitely many gauge equivalence classes of 
    minimizers.

{\bf b)} Let $\epsilon>0$.
         Then Condition C holds for any Palais-Smale sequence $(A_k,\Phi_k)$
         for $\YMH$ with $\YMH[A_k,\Phi_k]\ge\epsilon$ for all $k$.

{\bf c)} The set of critical values of $\YMH$ 
         is unbounded.

}

As was first pointed out by Gritsch [9],
there is a technical error in [14] on page 346, lines 11--13.
Here the Uhlenbeck compactness theorem is used to construct
a subsequence $A_k$ and a sequence of gauge transformations $g_k\in\eusG$ 
such that $g_k.A_k$ converges in $L^{1,2}$ 
to a Yang-Mills connection $A_\infty$.
Then the invariant Sobolev embedding $(L^{1,2})^H\to L^6$
is applied to the sequence $g_k.A_k$.
However, this embedding theorem can only be applied if all
the connections $g_k.A_k$ are known to be invariant under the
same action.
At this point in the proof,
the connection $g_k.A_k$ is only known to be invariant under
the gauge transformed action $g_k.\sigma_i$.
In the present paper we overcome this difficulty by working in 
the Morrey space $L^{1,2}_3$.
The advantage of using $L^{1,2}_3$ instead of the invariant Sobolev space
$(L^{1,2})^H$ as in [14] is that $L^{1,2}_3$ gives a topology on $\eusA$ 
while $(L^{1,2})^H$ only gives a topology on the subspace $\eusA^H_i$.
This is crucial in \S2.

It seems likely that one could work
directly in $L^{1,2}\cap L^6$ on manifolds of dimension $\le5$ or
in $L^{1,2}\cap L^6\cap C^{-1,1/2}$ on manifolds of any dimension
without introducing Morrey spaces.
One would then have to establish
first an embedding $(L^2)^H\to L^2\cap L^{-1,6}$
or $(L^2)^H\to L^2\cap L^{-1,6}\cap C^{-2,1/2}$
and then
an $L^{1,2}\cap L^6$ or $L^{1,2}\cap L^6\cap C^{-1,1/2}$
estimate for connections with $L^2\cap L^{-1,6}$ or
$L^2\cap L^{-1,6}\cap C^{-2,1/2}$ bounds respectively on curvature.

{\it Acknowledgements.} The author wishes to thank Gil Bor, Magnus Fontes,
Ursula Gritsch, Arne Meurman, Jaak Peetre, Lorenzo Sadun, and Karen Uhlenbeck 
for interesting discussions.
He also wishes to thank
the Department of Mathematics at University of Texas at Austin,
where some of this research was carried out, and in particular
Dan Freed for their hospitality.


\heading \S2. Proofs of Theorem~1 and Theorem~2

\noindent
We first review the basic setup for invariant gauge theory
as in [3], [6], [7], and [14].
The extended gauge group $\eusGhat$ is defined
as the group of ordered pairs $(\phi,h)$
such that $h\in H$ and $\phi:P\to P$ is a bundle map
that covers the isometry $\rho(h):X\to X$.
There is a short exact sequence
$$1\longrightarrow\eusG\longrightarrow\eusGhat
   \longrightarrow H\longrightarrow 1.
  \eqno(2.1)
$$
By a splitting of this short exact sequence we mean a smooth
left inverse $H\to\eusGhat$ of the homomorphism $\eusGhat\to H$.
Thus a splitting of (2.1) is the same as a lift $\sigma$ of $\rho$.

There are compatible group actions
$$\eqalign{
    \eusG\times\eusA&\to\eusA\cr
    \eusGhat\times\eusA&\to\eusA\cr
    H\times\eusA/\eusG&\to\eusA/\eusG.\cr
          }
$$
The first two actions are defined the obvious way.
The third action was introduced by Furuta [7] and is defined as follows:
Given $h\in H$ and $[A]\in\eusA/\eusG$,
let $\phi$ be any bundle map that covers the action of $h$ on $X$.
We then define $h.[A]=[\phi.A]\in\eusA/\eusG$.
As $A$ and $\phi$ are unique up to gauge transformations,
$\phi.A$ is unique up to gauge transformations.
Hence $[\phi.A]$ is a well-defined element of $\eusA/\eusG$.
We denote the fixed point set of the action 
$H\times\eusA/\eusG\to\eusA/\eusG$
by $(\eusA/\eusG)^H$.

For any $A\in\eusA$ these actions have isotropy subgroups
$$\eqalign{
    \Gamma_A & = \bigl\{g\in\eusG \,|\, g.A=A\bigr\} \cr
    \widehat\Gamma_A & = \bigl\{(\phi,h)\in\eusGhat \,|\, \phi.A=A\bigr\} \cr
    \Gamma_{[A]} & = \bigl\{h\in H  \,|\, h.[A]=[A]\bigr\} . \cr
          }
$$
These are compact Lie groups,
and they form a short exact sequence
$$1\longrightarrow\Gamma_A\longrightarrow\widehat\Gamma_A
   \longrightarrow \Gamma_{[A]}\longrightarrow 1.
   \eqno(2.2)
$$

We will work in Morrey spaces $L^{k,p}_d$.
These are function spaces well adapted for handling invariant problems.
They are discussed in some detail in Appendix~A;
here we briefly state their main properties.

\item $\bullet$

Functions in $L^{k,p}(X)$ that are invariant under a Lie group action on $X$,
all of whose orbits have codimension $\le d$, lie in $L^{k,p}_d(X)$.
(Lemma~A.1)

\item $\bullet$ 

The spaces $L^{k,p}_d(X)$ for $X$ of any dimension
satisfy ``the same'' multiplications, embeddings, and compact embeddings as
the space $L^{k,p}(X)$ for $X$ of dimension $d$.
(Lemma~A.2, Lemma~A.3, Lemma~A.4)

\item $\bullet$ 

Elliptic estimates hold in $L^{k,p}_d(X)$.
(Lemma~A.5, Remark~A.6)

\noindent
Using these properties it follows, in the usual way,
that for any $A_0\in\eusA$ there exists $\epsilon>0$ such that
$$\eusO_{A_0,\epsilon}
  = \bigl\{ A_0+a\in\eusA
    \,\big|\, \hbox{
          $d^*_{A_0}a=0$ and
          $\|a\|_{L^{1,2}_{3,A_0}(X,T^*X\otimes \Ad P)} <\epsilon$}
    \bigr\}
$$
is a local slice through $A_0$ for the action of $\eusG$ on $\eusA$.
More precisely, this means that the natural map
$$\eusG\times_{\Gamma_{A_0}}\eusO_{A_0,\epsilon}\to\eusA
  \eqno(2.3)
$$
is a $\eusG$-equivariant diffeomorphism onto its image.
In particular
$\eusO_{A_0,\epsilon}/\Gamma_{A_0}\to\eusA/\eusG$
is a homeomorphism onto its image.

Next we prove that the image of $\pi_i$ is closed.
This was established by Furuta [7], see also [3],
in the case when $G$ has trivial center 
and all connections in $\eusA^H_i$ are irreducible.
Here we extend his result to arbitrary compact $G$ and reducible connections.

\proclaim Lemma~2.1.

For any $A_0\in\eusA$ there exists $\epsilon>0$ 
such that if $A\in\eusO_{A_0,\epsilon}$, $\sigma$~is a lift of $\rho$,
and $A$ is invariant under $\sigma$,
then $A_0$ is invariant under $\sigma$.

\demo Proof.

If a connection $A$ is invariant under any lift $\sigma$ of $\rho$,
then $[A]\in(\eusA/\eusG)^H$.
The fixed point set $(\eusA/\eusG)^H$ is a closed subset of $\eusA/\eusG$
in the $L^{1,2}_3$ topology.
So if $[A_0]\notin(\eusA/\eusG)^H$,
then there exists $\epsilon>0$ such that
no connection in $\eusO_{A_0,\epsilon}$ is invariant under any lift $\sigma$
of $\rho$.
Thus we may assume that $[A_0]\in(\eusA/\eusG)^H$.
Then $\Gamma_{[A_0]}=H$, so the short exact sequence (2.2) becomes
$$1\longrightarrow\Gamma_A\longrightarrow\widehat\Gamma_A
   \longrightarrow H \longrightarrow 1.
   \eqno(2.4)
$$
It is straightforward to verify
that the action $\eusGhat\times\eusA\to\eusA$ restricts to an action
$\widehat\Gamma_{A_0}\times\eusO_{A_0,\epsilon}\to\eusO_{A_0,\epsilon}$.
Thus we can form the Banach manifold
$\widehat\eusG\times_{\widehat\Gamma_{A_0}}\eusO_{A_0,\epsilon}$.
It follows from (2.1) and (2.4) that the natural map
$\eusG\times_{\Gamma_{A_0}}\eusO_{A_0,\epsilon}
  \to 
  \widehat\eusG\times_{\widehat\Gamma_{A_0}}\eusO_{A_0,\epsilon}
$
is a diffeomorphism.
It then follows from (2.3) that the natural map
$$\widehat\eusG\times_{\widehat\Gamma_{A_0}}\eusO_{A_0,\epsilon}\to\eusA
  \eqno(2.5)
$$
is a diffeomorphism onto its image.
To prove the lemma,
it suffices to show that $\widehat\Gamma_A\subseteq \widehat\Gamma_{A_0}$.
Let $(\phi,h)\in\widehat\Gamma_A$.
Then $((1,1),A)$ and $((\phi,h),A)$ have the same image $A$ 
under the map (2.5).
Hence $(\phi,h)\in\widehat\Gamma_{A_0}$.
\enddemo

\proclaim Proposition~2.2.

$\pi_i(\eusA^H_i/\eusG^H_i)$ is a closed subset of $\eusA/\eusG$ in the $L^{1,2}_3$
topology.

\demo Proof.

Let $A_k\in\eusA^H_i$.
Assume that $[A_k]\to[A_\infty]$ in $L^{1,2}_3$ for some $A_\infty\in\eusA$.
This means that there exist $g_k\in\eusG$ such that $g_k.A_k\to A_\infty$
in $L^{1,2}_3$.
We have to show that $A_\infty$ is gauge equivalent 
to a connection in $\eusA^H_i$.

Let $\epsilon$ be as in Lemma~2.1.
By the slice theorem, we may assume that 
$g_k.A_k\in\eusO_{A_\infty,\epsilon}$ for large $k$.
As $A_k$ is invariant under $\sigma_i$,
we have that $g_k.A_k$ is invariant under $g_k.\sigma_i$.
It then follows from Lemma~2.1 that
$A_\infty$ is  invariant under $g_k.\sigma_i$ for large $k$.
Hence $g_k^{-1}.A_\infty$ is invariant under $\sigma_i$,
i.e.~$g_k^{-1}.A_\infty\in\eusA^H_i$.
\enddemo

\demo Proof of Theorem~1.

Let $A_k$ be a Palais-Smale sequence in $\eusA^H_i$.
It follows from Lemma~A.1 that
$$\|F_{A_k}\|_{L^2_3(X,\Lambda^2T^*X\otimes\Ad P)}
  \le M' \rlap{\qquad for all $k$}
$$
for some $M'>0$ and
$$\big\|d^*_{A_k}F_{A_k}\big\|_{L^{-1,2}_{3,A_k}(X,T^*X\otimes\Ad P)}
  \to 0 \rlap{\qquad as $k\to\infty$.}
$$
Arguing exactly as in the proof of Condition~C for Yang-Mills in three
dimensions, [17], [18] Proposition~4.5, [20], [23],
but using the Morrey space theory in Appendix~A instead of standard
Sobolev space theory,
we see that there exists a subsequence, which we also denote $A_k$,
and a sequence $g_k'\in\eusG$ such that
$g_k'.A_k$ converges in $L^{1,2}_3$ to a smooth
Yang-Mills connection $A'_\infty\in\eusA$.
It follows from Proposition~2.2 that there exists $g\in\eusG$ such that
$A_\infty=g.A'_\infty\in\eusA^H_i$.
Weak Condition~C follows with $g_k=g g'_k$.
\enddemo

\demo Proof of Theorem~2, sufficiency.

Assume that every Yang-Mills connection in $\eusA/\eusG$
has only finitely many preimages in $\eusA^H_i/\eusG^H_i$.
Let $A_k\in\eusA^H_i$ be a Palais-Smale sequence.
After passing to a subsequence there exists by Theorem~1
a sequence $g'_k\in\eusG$ such that $g'_k.A_k$ converges
in $L^{1,2}_3$ to a Yang-Mills connection $A'_\infty$.
Arguing as in the proof of Proposition~2.2
we may assume that
$${g'_k}^{-1}.A'_\infty\in\eusA^H_i$$
for large $k$.
In particular $[{g'_k}^{-1}.A'_\infty,\sigma_i]\in\pi_i^{-1}[A'_\infty]$
for large $k$.
By our assumptions, $\pi_i^{-1}[A'_\infty]$ is finite.
After passing to a subsequence we may therefore assume that
all $[{g'_k}^{-1}.A'_\infty,\sigma_i]$ are equal.
This means that for any integer $k_0$ there exists a sequence
$g_k\in\eusG^H_i$
such that
$g_k {g'_k}^{-1}.A'_\infty={g'_{k_0}}^{\!\!\!-1}.A'_\infty$ for all $k$.
Then
$$s_k=g'_{k_0}g_k {g'_k}^{-1}\in\Gamma_{A'_\infty}.$$
The group $\Gamma_{A'_\infty}$ is compact.
Hence, after passing to a subsequence,
$s_k$ converges in $C^\infty$ to some $s_\infty\in\Gamma_{A'_\infty}$
as $k\to\infty$.
It follows that
$$g_k.A_k = {g'_{k_0}}^{\!\!\!-1} s_k g'_k.A_k 
  \;\to\;
  {g'_{k_0}}^{\!\!\!-1} s_\infty.A'_\infty
  = {g'_{k_0}}^{\!\!\!-1}.A'_\infty\in\eusA^H_i
$$
in $L^{1,2}_3$ as $k\to\infty$.
Strong Condition~C follows with
$A_\infty={g'_{k_0}}^{\!\!\!-1}.A'_\infty$.

\demo Necessity.

If a Yang-Mills connection in $\eusA/\eusG$
has infinitely many preimages in $\eusA^H_i/\eusG^H_i$,
then these preimages form a Palais-Smale sequence
and should provide a counterexample to strong Condition~C.
To prove this, we have to verify that for any $[A]\in\eusA/\eusG$,
the set $\pi_i^{-1}([A])$ is a discrete subset of $\eusA^H_i/\eusG^H_i$.
This means that for any $A_0\in\eusA^H_i$
there exists $\epsilon>0$ such that if $A\in\eusA^H_i$,
$\|A-A_0\|_{L^{1,2}_{3,A_0}}\le\epsilon$, and $[A]=[A_0]\in\eusA/\eusG$,
then $[A,\sigma_i]=[A_0,\sigma_i]\in\eusA^H_i/\eusG^H_i$.
Consider such an $A$.
If $\epsilon$ is small enough, then it follows from the slice theorem 
that $A=\exp \psi.A_0$
for a unique
$\psi$ in $\Omega^0(\Ad P)$ that is $L^2$-perpendicular
to the null space of $d_{A_0}$  and has small $L^{2,2}_{3,A_0}$ norm.
By uniqueness, $\psi$ is invariant under $\sigma_i$,
so $\exp\psi\in\eusG^H_i$.
It follows that $[A,\sigma_i]=[A_0,\sigma_i]$.
\enddemo


\heading \S3. Proofs of Theorem~3 and Theorem~4

\noindent
In order to apply Theorem~2 we need a way of determining
$\pi^{-1}_i([A])$.
That is our next topic.
It follows from the definitions that
$g.\sigma_i=\sigma_i$ if and only if $g\in\eusG^H_i$.
In other words, $\sigma_i$ is invariant under $g$
if and only if $g$ is invariant under $\sigma_i$.
Using this one sees that there is a natural map
$$\bigcup_{i\in I}\eusA^H_i/\eusG^H_i 
  \quad\longrightarrow\quad
  \bigg\{\> (A,\sigma) \>\bigg|\>  \vcenter{
       \hbox{$A\in\eusA$, $\sigma$ is a lift of $\rho$, and}
       \hbox{$A$ is invariant under $\sigma$}
                                        } \>\bigg\}
  \bigg/ \eusG 
  \eqno(3.1)
$$
given by $\eusA^H_i \ni A\mapsto (A,\sigma_i)$,
and that this map is a 1--1 correspondence.
Here $\eusG$ acts by $g.(A,\sigma)=(g.A,g.\sigma)$.
Thus we can view an element of $\eusA^H_i/\eusG^H_i$
as the gauge equivalence class of a pair $(A,\sigma_i)$.
This is the reason we use the notation $[A,\sigma_i]$
for elements of $\eusA^H_i/\eusG^H_i$.

It follows from (3.1) that for any $A\in\eusA$ 
there is a 1--1 correspondence
$$\bigcup_{i\in I}\pi_i^{-1} [A]
  \quad\longrightarrow\quad
  \bigg\{\> \sigma \>\bigg|\>  \vcenter{
       \hbox{$\sigma$ is a lift of $\rho$ such that}
       \hbox{$A$ is invariant under $\sigma$}
                                        } \>\bigg\}
  \bigg/ \Gamma_A .
  \eqno(3.2)
$$
This can be reformulated in Lie theoretic terms.
Recall that a lift $\sigma$ of $\rho$ is the same as a splitting of 
the short exact sequence (2.1).
If there exists any lift $\sigma$ of $\rho$ 
such that $A$ is invariant under $\sigma$,
then $A\in(\eusA/\eusG)^H$, or equivalently $\Gamma_{[A]}=H$.
Then the short exact sequence (2.2) becomes
$$1\longrightarrow\Gamma_A\longrightarrow\widehat\Gamma_A
   \longrightarrow H\longrightarrow 1,
  \eqno(3.3)
$$
and a lift $\sigma$ of $\rho$ 
such that $A$ is invariant under $\sigma$ is the same as a splitting of 
this short exact sequence.
Combining this with (3.2) we get the following:

\proclaim Proposition~3.1.

Let $A\in\eusA$ with $[A]\in(\eusA/\eusG)^H$.
Then there is a 1--1 correspondence between the disjoint union
$\bigcup_{i\in I}\pi_i^{-1}[A]$ and the set of\/ $\Gamma_A$-conjugacy
classes of splittings of the short exact sequence (3.3).

{\it Example.}\enspace
For the trivial connection $A_0$ in Example~2 in \S1
the sequence (3.3) is
$$1\longrightarrow \SU(2) \longrightarrow \SU(2)\times\Spin(4)
   \longrightarrow \Spin(4) \longrightarrow 1 .
$$
Now $\Spin(4)$ is isomorphic to $\SU(2)\times\SU(2)$.
Thus this short exact sequence has three conjugacy classes of splittings,
given by
$\tau(x,y) = (1,x,y)$, $\omega_+(x,y) = (x,x,y)$ and $\omega_-(x,y)=(y,x,y)$.
It then follows from Proposition~3.1 that 
$[A_0]\in\eusA/\eusG$ has three preimages, 
$[A_0,\tau]$, $[A_0,\omega_+]$, and $[A_0,\omega_-]$,
in $\eusA^{\Spin(4)}_0/\eusG^{\Spin(4)}_0 
 \cup \eusA^{\Spin(4)}_1/\eusG^{\Spin(4)}_1
$.
In the notation of \S1 Example 2,
$\tau$ is gauge equivalent to $\sigma_0$ 
and $\omega_\pm$ are gauge equivalent to $\sigma_1$ .
$[A_0,\sigma_0]$ is the unique element in 
$\eusA^{\Spin(4)}_0/\eusG^{\Spin(4)}_0$
and $[A_0,\omega_\pm]$ are the points $t=\pm1$ in
$\eusA^{\Spin(4)}_1/\eusG^{\Spin(4)}_1$.

\demo Proof of Theorem~3.

If $A$ is irreducible, then $\Gamma_A=Z(G)$.
If $G$ is semisimple, then $Z(G)$ is finite.
It is easy to see that if $\Gamma_{A}$ is finite, 
then (3.3) has only finitely many splittings.
Theorem~3 now follows from Theorem~2 and Proposition~3.1.
\enddemo

\demo Proof of Theorem~4.

By Theorem~2 and Proposition~3.1, we only have to show that the
short exact sequence (3.3) has only finitely many $\Gamma_A$-conjugacy 
classes of splittings.
We denote the homomorphism $\widehat\Gamma_A\to H$ by $\mu$.
By Lemma~B.1 in Appendix~B
there are only finitely many $\widehat\Gamma_A$-conjugacy classes of 
homomorphisms $H\to\widehat\Gamma_A$.
Thus we have to show that any $\widehat\Gamma_A$-conjugacy class 
of homomorphisms $H\to\widehat\Gamma_A$ contains
only finitely many $\Gamma_A$-conjugacy classes of splittings of (3.3).

If $\sigma\in\Hom(H,\widehat\Gamma_A)$ is a splitting,
then the other splittings in the $\widehat\Gamma_A$-conjugacy class of $\sigma$
are of the form $k\sigma k^{-1}$
with $k\in\mu^{-1}(Z(H))$.
Two such splittings $k_1\sigma k_1^{-1}$ and $k_2\sigma k_2^{-1}$ 
lie in the same $\Gamma_A$-conjugacy class 
if $k_1 k_2^{-1}\in \Gamma_A\,Z(\widehat\Gamma_A)$.
Thus it suffices to show that $\Gamma_A\,Z(\widehat\Gamma_A)$ 
is a subgroup of $\mu^{-1}(Z(H))$ of finite index.
It is a closed normal subgroup. We have
$\mu^{-1}(Z(H))/\Gamma_A\,Z(\widehat\Gamma_A) 
 = Z(H)/\mu (Z(\widehat\Gamma_A))$.
Next we analyze the Lie algebras of the groups involved.

The short exact sequence (3.3) of Lie groups induces a short exact sequence
$$0\longrightarrow\gamma_A\longrightarrow\widehat\gamma_A
   \mapright{\mu_*}\eufh\longrightarrow 0
$$
of Lie algebras.
Recall that if $G$ is a compact Lie group,
then its Lie algebra $\eufg$ has a unique decomposition 
as a direct sum $\eufz(\eufg)\oplus\eufg_\ss$
of an abelian Lie algebra and a semisimple Lie algebra;
the abelian part is the center of $\eufg$.
An ideal or quotient of an abelian or semisimple Lie algebra
is abelian or semisimple respectively.
Hence we have short exact sequences
$$0\longrightarrow\eufz(\gamma_A)\longrightarrow\eufz(\widehat\gamma_A)
   \mapright{\mu_*}\eufz(\eufh)\longrightarrow 0
$$
and
$$0\longrightarrow(\gamma_A)_\ss\longrightarrow(\widehat\gamma_A)_\ss
   \mapright{\mu_*}\eufh_\ss\longrightarrow 0
$$
of Lie algebras.
It follows that $Z(H)$ and $\mu(Z(\widehat\Gamma_A))$
have the same Lie algebra $\eufz(\eufh)=\mu_*\eufz(\widehat\gamma_A)$.
Hence their quotient $Z(H)/\mu(Z(\widehat\Gamma_A))$ is a discrete group.
But it is also compact, so it is finite.
The theorem follows.
\enddemo


\heading \S4. Proof of Corollary~2

\noindent
Let $\rho$ be the $\Spin(3)$ action on $S^2\times S^2$ 
given by the standard action on each factor.
Corollary~2 is established by applying Theorem~4 to
spaces of connections that are invariant under lifts $\sigma$ of $\rho$,
using the ``the mountain pass lemma'' argument for
paths between $\U(1)$-reducible Yang-Mills connections.
To carry out this we have to first classify the lifts $\sigma$ of $\rho$
to principal $\SU(2)$-bundles over $S^2\times S^2$,
then classify the $\U(1)$-reducible Yang-Mills connections
that are invariant under these lifts,
and finally determine the index and nullity
of the Hessian of the Yang-Mills functional 
on the normal bundles of the gauge orbits of these connections.
This is done in Lemma~4.1, Lemma~4.2, and Lemma~4.3 below.

The principal $\SU(2)$ bundles over $S^2\times S^2$ are classified by the
second Chern number.
We let $P_d$ denote the bundle with second Chern number $d$.
We let $\eusA_d$ and $\eusG_d$ denote the space of connections
and the group of gauge transformations on $P_d$.

\proclaim Lemma~4.1.

The lifts $\sigma$ of the action $\rho$ of $\Spin(3)$ on $S^2\times S^2$
to $P_d$ are classified, up to gauge equivalence,
by pairs $(w_+,w_-)$ of nonnegative integers
such that $w_+\equiv w_-$ {\rm(mod 2)} and $d=(w_+^2-w_-^2)/2$.

We denote the corresponding lifts by $\sigma_{w_+,w_-}$
and we denote the space of connections and the group of gauge transformations
that are invariant under $\sigma_{w_+,w_-}$ by $\eusA^{\Spin(3)}_{w_+,w_-}$
and $\eusG^{\Spin(3)}_{w_+,w_-}$.

\demo Proof.

This follows from [21] \S3,
except the formula for $d$,
which follows by examining the examples $[A_{p,q},\tau_{p,q}]$ below.
\enddemo

The numbers $w_+$ and $w_-$ have the following geometric meaning.
For $(x,y)\in S^2\times S^2$, let $\theta$ denote the angle between
$x$ and $y$, viewed as unit vectors in $\IR^3$.
Then $\theta\in[0,\pi]$ parametrizes the space of orbits of the action $\rho$.
For $\theta\in(0,\pi)$, the orbit is a copy of $\SO(3)=\Spin(3)/\IZ_2$.
For $\theta=0$ and $\pi$, the orbit is a copy of $S^2=\Spin(3)/\Spin(2)$.
The orbit $\theta=0$ consists of all pairs $(x,x)$ 
and the orbit $\theta=\pi$ of all pairs $(x,-x)$ with $x\in S^2$.
We denote these two orbits by $S^2_+$ and $S^2_-$.
Each point in the orbits $S^2_+$ and $S^2_-$ has isotropy subgroup $\Spin(2)$.
Thus $\Spin(2)$ acts on the fibers over these two orbits.
These actions give two homomorphisms from $\Spin(2)$ to $\SU(2)$,
well-defined up to conjugacy.
Such conjugacy classes of homomorphisms are classified
by nonnegative integral weights. This gives the weights $w_+$ and $w_-$.
For more details, see [21] \S3.

Next we construct the $\U(1)$-reducible Yang-Mills connections in
$\eusA^{\Spin(3)}_{w_+,w_-}/\eusG^{\Spin(3)}_{w_+,w_-}$.
Principal $\U(1)$-bundles over $S^2$ are classified by the first Chern number.
The standard $\Spin(3)$-action on $S^2$ has unique lifts to these bundles
up to gauge equivalence.
There are unique invariant Yang-Mills connections on these bundles.
We take the tensor product of the Chern number $p$ and $q$ principal
$\U(1)$-bundles over the first and second factor of $S^2\times S^2$
and we then form the associated principal $\SU(2)$-bundle.
It has second Chern number $2pq$.
The construction gives an $\Spin(3)$-action $\tau_{p,q}$ on the bundle
that covers $\rho$
and a $\U(1)$-reducible invariant Yang-Mills connection $A_{p,q}$.
The action $\tau_{p,q}$ has $w_+=|p+q|$ and $w_-=|p-q|$.
By Lemma~4.1, $\tau_{p,q}$ is thus gauge equivalent to
$\sigma_{|p+q|,|p-q|}$.
Hence  $[A_{p,q},\tau_{p,q}]$ defines an element in
$\eusA^{\Spin(3)}_{|p+q|,|p-q|}/\eusG^{\Spin(3)}_{|p+q|,|p-q|}$.
We have $[A_{p,q},\tau_{p,q}]=[A_{-p,-q},\tau_{-p,-q}]$;
otherwise these elements are distinct.

Another lift $\omega$ of $\rho$ to the trivial bundle 
$S^2\times S^2\times \SU(2) \to S^2\times S^2$ is obtained by 
letting $\Spin(3)$ act in the standard way on the first two factors
and by left multiplication on the third factor.
This action has $w_+=w_-=1$.
By Lemma~4.1, $\omega$ is thus gauge equivalent to $\sigma_{1,1}$.
The trivial connection $A_{0,0}$ is invariant under $\omega$.
Thus $[A_{0,0},\omega]$ defines an element of
$\eusA^{\Spin(3)}_{1,1}/\eusG^{\Spin(3)}_{1,1}$.
The following lemma says that there are no other reducible invariant
Yang-Mills connections.

\proclaim Lemma~4.2.

{
The following are the elements of
$\eusA^{\Spin(3)}_{w_+,w_-}/\eusG^{\Spin(3)}_{w_+,w_-}$ 
that are given by\/ $\U(1)$-reducible Yang-Mills connections:

$w_+,w_->0$, $(w_+,w_-)\ne(1,1)${\rm:}\quad
$[A_{p,q},\tau_{p,q}]$ and $[A_{q,p},\tau_{q,p}]$
with $p=(w_++w_-)/2$ and $q=(w_+-w_-)/2$,

$(w_+,w_-)=(1,1)${\rm:}\quad
$[A_{1,0},\tau_{1,0}]$, $[A_{0,1},\tau_{0,1}]$, and $[A_{0,0},\omega]$,

$w_+>0$, $w_-=0${\rm:}\quad
$[A_{p,p},\tau_{p,p}]$ with $p=w_+/2$,

$w_+=0$, $w_->0${\rm:}\quad
$[A_{p,-p},\tau_{p,-p}]$ with $p=w_-/2$,

$w_+=w_-=0${\rm:}\quad
$[A_{0,0},\tau_{0,0}]$.

}

\demo Proof.

The reducible Yang-Mills connections in $\eusA_d/\eusG_d$ are
$[A_{p,q}]$ with $2pq=d$.
To establish the lemma we only have to find the preimages
in $\bigcup_{(w_+,w_-)}  \eusA^{\Spin(3)}_{w_+,w_-}/\eusG^{\Spin(3)}_{w_+,w_-}
   $
of $[A_{p,q}]$.
By Proposition~3.1 these are classified by 
the conjugacy classes of splittings
of the sequence (3.3) for $A_{p,q}$.

For $[A_{p,q}]$ with $(p,q)\ne 0$, the sequence is
$$1\to \U(1) \to \U(1)\times\Spin(3) \to \Spin(3) \to 1.$$
This sequence has exactly one splitting.
The corresponding preimage is $[A_{p,q},\tau_{p,q}]$.

For $[A_{0,0}]$ the sequence is 
$$1\to \SU(2) \to \SU(2)\times\Spin(3) \to \Spin(3) \to 1.$$
This sequence has exactly two conjugacy classes of splittings.
The corresponding preimages are
$[A_{0,0},\tau_{0,0}]$ and $[A_{0,0},\omega]$.
\enddemo

In the following we view elements of 
$\eusA^{\Spin(3)}_{w_+,w_-}/\eusG^{\Spin(3)}_{w_+,w_-}$
as $\eusG^{\Spin(3)}_{w_+,w_-}$-orbits in
$\eusA^{\Spin(3)}_{w_+,w_-}$.

\proclaim Lemma~4.3.

The index and nullity of the Hessian of the Yang-Mills functional
on the normal bundle in $\eusA^{\Spin(3)}_{w_+,w_-}$
of the $\eusG^{\Spin(3)}_{w_+,w_-}$-orbit
of a $\U(1)$-reducible Yang-Mills connection
in $\eusA^{\Spin(3)}_{w_+,w_-}$ is as follows:
$$\vbox{\halign{
  $#$\hfil &
    & \quad \hfil$#$\hfil \cr
  [A,\sigma] & \hbox{\it index} & \hbox{\it nullity} \cr
  \noalign{\smallskip}
  [A_{0,\pm2} , \tau_{0,\pm2}] & 0 & 2 \cr
  [A_{\pm2,0} , \tau_{\pm2,0}] & 0 & 2 \cr
  [A_{p,p\pm1} , \tau_{p,p\pm1}] & 2 & 0 \cr
  [A_{p,q} , \tau_{p,q}] & 0 & 0 & \rlap{\it for other $(p,q)$} \cr
  [A_{0,0} , \omega] & 0 & 0\rlap{.} \cr
         }}
$$

\demo Proof.

The Hessian of the Yang-Mills functional is
$H_A a = d_A^* d_A a + *[*F_A,a]$.
The spectral decomposition of the operator
$(H_A+d_Ad_A^*)a = \Delta_A a + *[*F_A,a]$
has been carried out in [22] Appendix~A.
There $T_A\eusA_d=\Omega^1(S^2\times S^2,\Ad P_d)$ is decomposed
as an orthogonal direct sum of six subspaces
that are invariant under $H_A+d_Ad_A^*$.
On these subspaces the operator is written in terms of
the scalar Laplacian on $S^2$ and $\overline\partial$-Laplacians on
holomorphic line bundles over $S^2$.
It is shown in [22] Appendix~A,
although never stated explicitly,
that the spectral decomposition of $H_A+d_Ad_A^*$ 
on $T_A\eusA_d$ is as follows:
$$\vbox{\halign{$#$\hfil\quad&&\quad$#\hfil$\quad\cr
          [A]& \rm eigenvalue & \rm multiplicity \cr
          \noalign{\smallskip}
          [A_{p,q}]
          & k^2+k+l^2+l
          & 2(2k+1)(2l+1)
          & k\ge 1,l\ge 0 \cr
          & k^2+k+l^2+l
          & 2(2k+1)(2l+1)
          & k\ge 0,l\ge 1 \cr
          & k^2 + k -p^2 + l^2 + l - q^2
          & 2(2k+1)(2l+1)
          & k\ge |p-1|, l\ge |q| \cr
          & k^2 + k - p^2 +l^2 + l - q^2 
          & 2(2k+1)(2l+1)
          & k\ge |p+1| ,l\ge |q| \cr
          & k^2 + k -p^2 + l^2 +l - q^2 
          & 2(2k+1)(2l+1)
          & k\ge|p|, l\ge|q-1| \cr
          & k^2 + k -p^2 + l^2 + l - q^2 
          & 2(2k+1)(2l+1)
          & k\ge |p|, l\ge |q+1|.  \cr
     }}
$$
Multiplicities for repeated eigenvalues should be added.
It is straightforward to keep track of our $\Spin(3)$-actions throughout 
[22] Appendix~A.
The action of $\Spin(3)$ on $\eusA_d$ induces actions on the eigenspaces
as follows:
$$\vbox{\halign{$#$\hfil\quad&&\quad$#\hfil$\quad\cr
          [A,\sigma] & \rm eigenvalue & \rm eigenspace \cr
          \noalign{\smallskip}
          [A_{p,q},\tau_{p,q}]
          & k^2+k+l^2+l
          & 2 V_k \otimes V_l
          & k\ge 1,l\ge 0 \cr
          & k^2+k+l^2+l
          & 2 V_k \otimes V_l
          & k\ge 0,l\ge 1 \cr
          & k^2 + k -p^2 + l^2 + l - q^2
          & 2 V_k\otimes V_l
          & k\ge |p-1|, l\ge |q| \cr
          & k^2 + k - p^2 +l^2 + l - q^2 
          & 2 V_k\otimes V_l
          & k\ge |p+1| ,l\ge |q| \cr
          & k^2 + k -p^2 + l^2 +l - q^2 
          & 2 V_k \otimes V_l
          & k\ge|p|, l\ge|q-1| \cr
          & k^2 + k -p^2 + l^2 + l - q^2 
          & 2 V_k \otimes V_l
          & k\ge |p|, l\ge |q+1|  \cr
          }}
$$
$$\vbox{\halign{$#$\hfil & \qquad$#$\hfil & \qquad\hfil$#$\hfil 
                & \qquad$#$\hfil \cr
          [A_{0,0},\omega]
          & k^2+k+l^2+l
          & 2 V_k \otimes V_l \otimes V_1
          & k\ge 1,l\ge 0 \cr
          & k^2+k+l^2+l
          & 2 V_k \otimes V_l \otimes V_1
          & k\ge 0,l\ge 1. \cr
        }}
$$
Here $V_k$ denotes the spin-$k$ representation of $\Spin(3)$.
The $\Spin(3)$-invariant subspaces of the above eigen\-spaces
give the spectral decomposition of $H_A+d_Ad_A^*$ on $T_A\eusA^{\Spin(3)}_{w_+,w_-}$:
$$\vbox{\halign{$#$\hfil & \quad$#$\hfil & \quad\hfil$#$\hfil 
                & \quad$#$\hfil \cr
          [A,\sigma] & \rm eigenvalue & \rm multiplicity \cr
          \noalign{\smallskip}
          [A_{p,q},\tau_{p,q}]\hbox{ with } |p|>|q| \qquad 
          & 2k^2 + 2k & 4 & k\ge 1 \cr
          & p^2 - 2|p| - q^2 & 2 \cr
          & p^2 +  2|p| - q^2 & 6 \cr
          & 2k^2 + 2k - p^2 - q^2 & 8 & k\ge |p|+1 \cr
          \noalign{\smallskip}
          [A_{p,q},\tau_{p,q}]\hbox{ with } |p|=|q|\ge1
          & 2k^2 + 2k & 4 & k\ge 1 \cr
          & 2|p| & 4 \cr
          & 2k^2 + 2k - 2p^2  & 8 & k\ge |p|+1 \cr
          \noalign{\smallskip}
          [A_{0,0},\tau_{0,0}]\qquad
          & 2k^2+2k & \qquad 12\qquad & k\ge 1 \cr
          \noalign{\smallskip}
          [A_{0,0},\omega] 
          & 2k^2 + 2k & 4 & k\ge 1 \cr
          & 2 & 4 \cr
          & 2k^2 & 8 & k\ge 2 . \cr
         }}
$$
The spectrum is unchanged if we replace $(p,q)$ by $(q,p)$.
Therefore we may assume that $|p|\ge|q|$.
Arguing as above, one sees that the spectral decomposition
of $\Delta_A$ on $\Omega^0(\Ad P_d)^{\Spin(3)}$ is as follows:
$$\vbox{\halign{$#$\hfil & \quad$#$\hfil & \quad\hfil$#$\hfil 
                & \quad$#$\hfil \cr
          [A,\sigma] & \rm eigenvalue & \rm multiplicity \cr
          \noalign{\smallskip}
          [A_{p,q},\tau_{p,q}]\hbox{ with $|p|\ge|q|$ and $|p|\ge 1$}
          & 2k^2+2k
          & 1
          & k\ge 0 \cr
          & 2k^2+2k-p^2-q^2 & 2
          & k\ge |p| \cr
          \noalign{\smallskip}
          [A_{0,0},\tau_{0,0}]\qquad
          & 2k^2+2k
          &  \qquad 3 \qquad
          & k\ge 0 \cr
          \noalign{\smallskip}
          [A_{0,0},\omega]
          & 2k^2 + 2k & 1 & k\ge1 \cr
          & 2k^2 & 2 & k\ge 1. \cr
         }}
$$
After removing the zero eigenspace
we have the spectral decomposition of $d_Ad_A^*$ on
 $T_A(\eusG^{\Spin(3)}_{w_+,w_-}.A)$.
Removing this from the spectral decomposition of $H_A+d_Ad_A^*$ on 
$T_A\eusA^{\Spin(3)}_{w_+,w_-}$ we finally get the spectral decomposition of
$H_A$ on $T_A^\perp(\eusG^{\Spin(3)}_{w_+,w_-}.A)$:
$$\vbox{\halign{$#$\hfil & \quad$#$\hfil & \quad\hfil$#$\hfil 
                & \quad$#$\hfil \cr
          [A,\sigma] & \rm eigenvalue & \rm multiplicity \cr
          \noalign{\smallskip}
          [A_{p,q},\tau_{p,q}]\hbox{ with }|p|>|q|
          & 2k^2 + 2k & 3 & k\ge 1 \cr
          & p^2 - 2|p| - q^2 & 2 \cr
          & p^2 +  2|p| - q^2 & 4 \cr
          &2k^2 + 2k - p^2 - q^2 & 6 & k\ge |p|+1 \cr
         }}
$$
$$\vbox{\halign{$#$\hfil & \qquad$#$\hfil & \qquad\hfil$#$\hfil 
                & \qquad$#$\hfil \cr
          [A_{p,q},\tau_{p,q}]\hbox{ with }|p|=|q|\ge1
          & 2k^2 + 2k & 3 & k\ge 1 \cr
          & 2|p| & 2 \cr
          & 2k^2 + 2k - 2p^2  & 6 & k\ge |p|+1 \cr
          \noalign{\smallskip}
          [A_{0,0},\tau_{0,0}]
          & 2k^2 + 2k & 9 & k\ge 1 \cr
          \noalign{\smallskip}
          [A_{0,0},\omega] 
          & 2k^2 + 2k & 3 & k\ge 1 \cr
          & 2 & 2 \cr
          & 2k^2 & 6 & k\ge 2 .\cr
         }}
$$
The only nonpositive eigenvalues are $-1$ in the case $|p|=|q|+1$
and 0 in the case $|p|=2$, $q=0$. Both have multiplicity 2.
The lemma follows.
\enddemo

\demo Proof of Corollary~2.

Let $d$ be an even integer $\ne\pm2,\pm4$.
Then there exist nonnegative integral weights $w_+$ and $w_-$,
as in Lemma~4.1, with $w_+,w_-\ne0,1$, $(w_+,w_-)\ne(2,2)$
and $d=(w_+^2-w_-^2)/2$.
By Lemma~4.2
there are exactly two $\eusG^{\Spin(3)}_{w_+,w_-}$-orbits
of $\U(1)$-reducible Yang-Mills connections in $\eusA^{\Spin(3)}_{w_+,w_-}$.
By Lemma~4.3 the Hessian of the Yang-Mills functional
is positive definite on the  normal bundles of these orbits.
It then follows from Theorem~4,
using ``the mountain pass lemma'' argument,
that there exists a Yang-Mills connection in $\eusA^{\Spin(3)}_{w_+,w_-}$
with energy greater than the energy of the two reducible ones.
This Yang-Mills connection has to be 
irreducible non-(anti-)selfdual and has second Chern number $d$.
\enddemo


\heading \S5. Proof of Corollary~3

\noindent
\demo Proof of Corollary~3.

The minimum of $\YMH$ is  0,
and the minimizers are the pairs $(A,\Phi)$ where $A$ has
covariantly constant curvature and $\Phi=*F_A$.
Such connections are $U(1)$-reducible;
the reduction is given by the spectral decomposition of $F_A$.
The $\U(1)$-bundles over $S^2$ are classified by the first Chern number.
On each such bundle there is a unique connection with covariantly 
constant curvature.
These connections on the $\U(1)$-bundles with first Chern numbers
$k$ and $-k$ give the same $\SU(2)$-connection.
Thus the minima of $\YMH$ are naturally indexed 
by nonnegative integers $k$.
This completes the proof of Part a).

The Hopf fibration gives a free action 
$\rho:\SO(2)\times S^3\to S^3$.
Here we take $S^3$ to be the 3-sphere of radius 2.
Then the quotient space is the 2-sphere of radius 1.
There is a unique principal $\SU(2)$-bundle $P$ over $S^3$.
The action $\rho$ lifts in a unique way, up to gauge equivalence,
to an action $\sigma:\SO(2)\times P\to P$.
Thus we omit the index $i$ and write $\eusA^{\SO(2)}$ and $\eusG^{\SO(2)}$.
It is a straightforward calculation to verify that,
apart from a constant,
$\YMH$ is the dimensional reduction
of $\YM$ under this action.
In particular, the space of pairs
$(A,\Phi)$ over $S^2$ can be identified with $\eusA^{\SO(2)}$
and the group of gauge transformations over $S^2$ can be identified with
$\eusG^{\SO(2)}$.
Thus Condition~C for $\YMH$ over $S^2$ is equivalent to strong
Condition~C for $\YM$ on $\eusA^{\SO(2)}$.
We will state the rest of the proof in terms of $\YM$ on
$\eusA^{\SO(2)}/\eusG^{\SO(2)}$.

We need a slight extension of the general theory.
The following theorem is proven exactly the same way as Theorem~2.

\proclaim Theorem~$2'\!.$

Let $a,b\in[0,\infty)$ with $a<b$.
Assume that all orbits of the action of $H$ on $X$ have codimension $\le 3$.
Then strong Condition~C holds for every
Palais-Smale sequence $A_k\in\eusA^H_i$, with $\YM[A_k]\in[a,b]$ for all $k$,
if and only if every Yang-Mills connection
$[A]\in\eusA/\eusG$, with $\YM[A]\in[a,b]$,
has only finitely many preimages in $\eusA^H_i/\eusG^H_i$
under the map $\pi_i$.

\penalty-100
Combining this with Proposition~3.1 we get the following.

\proclaim Theorem~$3'.\!$

Let $a,b\in[0,\infty)$ with $a<b$.
If all orbits of the action of $H$ on $X$ have codimension $\le 3$,
$G$ is semisimple,
and every Yang-Mills connection $A\in\eusA^H_i$ 
with $\YM[A]\in[a,b]$ is irreducible,
then strong Condition~C holds for any Palais-Smale sequence 
$A_k\in\eusA^H_i$ with $\YM[A_k]\in[a,b]$ for all $k$.

The trivial connection is the only reducible 
$\SU(2)$ Yang-Mills connection over $S^3$.
Part b) thus follows from Theorem~$3'$.

Finally we sketch the proof of Part c).
Let $\eusG_0$ be the subgroup of $\eusG$ consisting of gauge transformations
that act trivially on the fiber over a fixed point in $S^3$.
Let $\eusG^{\SO(2)}_0$ be the subgroup of $\eusG^{\SO(2)}$ consisting 
of gauge transformations that act trivially on the fibers over a fixed
$\SO(2)$-orbit in $S^3$.
Then $\eusA/\eusG_0$ and $\eusA^{\SO(2)}/\eusG^{\SO(2)}_0$ 
are Hilbert manifolds.
Using the above form of strong Condition~C
one can then establish a Morse theory
for a small perturbation of $\YM$ on 
$\YM^{-1}([\epsilon,\infty))\subseteq\eusA^{\SO(2)}/\eusG^{\SO(2)}_0$
for any $\epsilon>0$.

On the other hand,
Condition~C holds for $\YM$ on $\eusA/\eusG_0$.
The trivial connection is a nondegenerate critical point of $\YM$
on $\eusA/\eusG_0$.
Hence there exists $\epsilon>0$ such that the trivial connection
is the only critical point of $\YM$ in 
$\YM^{-1}([0,\epsilon])\subseteq\eusA/\eusG_0$.
Thus the preimages of the trivial connection are the only critical points
of $\YM$  in 
$\YM^{-1}([0,\epsilon])\subseteq\eusA^{\SO(2)}/\eusG^{\SO(2)}_0$.
These are also nondegenerate;
it is not hard to determine the spectral decompositions of the Hessians.
It follows from Part a) that the trivial connection in $\eusA/\eusG_0$
has countably infinitely many preimages in $\eusA^{\SO(2)}/\eusG^{\SO(2)}_0$.
This also follows from Proposition~3.1;
for the trivial connection the sequence (3.3) is
$$1\to \SU(2)\to\SU(2)\times\SO(2)\to\SO(2)\to1.$$
We conclude that 
$\YM^{-1}([0,\epsilon])\subseteq \eusA^{\SO(2)}/\eusG^{\SO(2)}_0$ 
is the disjoint union of countably infinitely many 0-cells.

This gives a decomposition
of $\eusA^{\SO(2)}/\eusG^{\SO(2)}_0$ as a CW-complex
with infinitely many 0-cells.
The Hilbert manifold $\eusA^{\SO(2)}/\eusG^{\SO(2)}_0$ is connected, 
so there then have to be infinitely many 1-cells.
These are given by critical points of index 1.
By the above form of strong Condition~C, 
they can not have uniformly bounded energy.
\enddemo


{\ninepoint

\heading \tenbf Appendix~A. Morrey spaces

\noindent
The Morrey space $L^p_\mi(\IR^n)$,
with $p\in[1,\infty)$ and $\mi\in[0,\infty)$,
is defined as the space of all $f\in L^p(\IR^n)$
such that
$$\sup_{\rho\in(0,1]}\sup_{x\in\IR^n} \rho^{\mi-n}
               \| f\|_{L^p(B_\rho(x))} ^p 
  <\infty .
$$
It is a Banach space with norm
$$\|f\|_{L^p_\mi(\IR^n)}^p
  = \|f\|_{L^p(\IR^n)}^p
          + \sup_{\rho\in(0,1]}\sup_{x\in\IR^n}
               \rho^{\mi-n} \|f\|_{L^p(B_\rho(x))} ^p .
  \eqno(A.1)
$$
The Morrey space $L^{k,p}_\mi(\IR^n)$,
with $k$ a positive integer, $p\in[1,\infty)$ and $\mi\in\IR$,
is defined as the space
of all $f\in L^{k,p}(\IR^n)$ such that
$\partial^\alpha_x f\in L^p_\mi(\IR^n)$ for all $\alpha$
with $|\alpha|\le k$.
It is a Banach space with norm
$$\|f\|_{L^{k,p}_\mi(\IR^n)}^p
  = \sum_{|\alpha|\le k} \|\partial^\alpha_x f\|_{L^p_\mi(\IR^n)}^p .
$$
The Morrey space $f\in L^{-k,p}_\mi(\IR^n)$,
with $k$ a positive integer, $p\in[1,\infty)$ and $\mi\in\IR$,
is defined as the space of all $f\in L^{-k,p}(\IR^n)$
for which there exist $g_\alpha\in L^p_\mi (\IR^n)$, $|\alpha|\le k$,
such that
$$f=\sum_{|\alpha|\le k} \partial^\alpha_x g_\alpha.
  \eqno(A.2)
$$
It is a Banach space with norm
$$\|f\|_{L^{-k,p}_\mi(\IR^n)}^p = \inf_{(g_\alpha)}\sum_{|\alpha|\le k}
  \|g_\alpha\|_{L^p_\mi(\IR^n)}^p
$$
where we take infimum over all collections $(g_\alpha)_{|\alpha|\le k}$
of functions in $L^p_\mi(\IR^n)$ that satisfy (A.2).
Note that $L^{k,p}_0(\IR^n) = L^{k,\infty}(\IR^n)$
and $L^{k,p}_\mi(\IR^n) = L^{k,p}(\IR^n)$ for $\mi\ge n$.

Let $\Omega$ be an open subset of $\IR^n$.
Then the Morrey space $L^{k,p}_\mi(\Omega)$,
with $k$ an integer, $p\in[1,\infty)$ and $\mi\in\IR$,
is defined as the space of all $f\in L^{k,p}(\Omega)$
such that there exists $g\in L^{k,p}_\mi(\IR^n)$
with $f=g\big|_\Omega$.
For negative $k$ the restriction is to be understood in the sense
of distributions.
It is a Banach space with norm
$$\|f\|_{L^{k,p}_\mi(\Omega)}
  = \inf_g  \|g\|_{L^{k,p}_\mi(\IR^n)}
  \eqno(A.3)
$$
where we take infimum over all $g\in L^{k,p}_\mi(\IR^n)$
such that $f=g\big|_\Omega$.

Let $X$ be an $n$-dimensional smooth compact Riemannian manifold.
Then the Morrey space $L^p_\mi(X)$,
with $k$ a positive integer and $\mi\in\IR$,
is defined as the space
of all $f\in L^p(X)$ such that the norm
$$\|f\|_{L^p_{\mi}(X)}^p
  =\|f\|_{L^p(X)}^p
   + \sup_{\rho\in(0,\rho_0/2]}\sup_{x\in X}
   \rho^{\mi-n} \|f\|_{L^p(B_\rho(x))}^p 
  \eqno(A.4)
$$
is finite;
here $\rho_0$ is the injectivity radius of $X$.
This norm is invariant under isometries of $X$.
Let $E$ be a Euclidean vector bundle over $X$.
Let $A$ be a smooth connection on $E$.
Then the Morrey space $L^{k,p}_\mi(X,E)$,
with $k$ a nonnegative integer, $p\in[1,\infty)$ and $\mi\in\IR$,
is defined as the space of all $s\in L^{k,p}(X,E)$
such that the norm
$$\|s\|_{L^{k,p}_{\mi,A}(X,E)}^p
  =\sum_{j=0}^k \bigl\|(|\nablaA^js|)\bigr\|_{L^p_\mi(X)}^p
  \eqno(A.5)
$$
is finite.
The Morrey space $L^{-k,p}_\mi(X,E)$,
with $k$ a positive integer, $p\in[1,\infty)$ and $\mi\in\IR$,
is defined as the space of all $s\in L^{-k,p}(X,E)$
such that the norm
$$\|s\|_{L^{-k,p}_{\mi,A}(X,E)}^p
  = \inf_{(t_j)}\sum_{j=0}^k
       \bigl\|(|t_j|)\bigr\|_{L^p_\mi(X)}^p
  \eqno(A.6)
$$
is finite;
here we take infimum over all $(t_0,\ldots,t_k)$ such that
$t_j\in L^p_\mi\bigl(X,(T^*X)^{\otimes j}\otimes E\bigr)$
and
$s = \sum_{j=0}^k (\nablaA^*)^jt_j$.
Different smooth connections $A$ give equivalent norms.

We can get equivalent norms for $L^{k,p}_\mi(X,E)$
as follows.
Choose an atlas of good local coordinate charts
$\Phi_\nu:\Omega_\nu\to \IR^n$
and good lifts $\Psi_\nu:E\big|_{\Omega_\nu}\to\IR^n\times\IR^N$.
Here good means that the maps can be extended smoothly to
neighborhoods of $\overline\Omega_\nu$.
Then
$$\|s\|^p =
  \sum_\nu \bigl\| (\Psi^{-1}_\nu)^* s
                \|_{L^{k,p}_\mi(\Phi_\nu(\Omega_\nu),\IR^N)}^p 
  \eqno(A.7)
$$
is a norm for $L^{k,p}_\mi(X,E)$,
with $k\in\IZ$, $p\in[1,\infty)$ and $\mi\in\IR$.
Using this norm 
analysis in Morrey spaces on manifolds can be reduced to analysis
in Morrey spaces on $\IR^n$.

\proclaim Lemma~A.1.

{
Let $H$ be a compact Lie group that acts smoothly
on $X$, in such a way that all $H$-orbits have codimension $\le d$,
and that acts smoothly on $E$ covering the action on $X$.
If $s\in L^{k,p}(X,E)$, with $k\in\IZ$ and $p\in[1,\infty)$,
is $H$-invariant, then $s\in L^{k,p}_d(X,E)$.
If $A$ is an $H$-invariant connection on $E$,
then
$$\|s\|_{L^{k,p}_{d,A}(X,E)} \le  c\|s\|_{L^{k,p}_A(X,)}.$$
The constant $c$ does not depend on $A$ and $s$.
}

We may assume that the connection $A$ in the norms (A.5) and (A.6)
is invariant under $H$.
Then the norms (A.4)--(A.6) are invariant under the action of $H$.

\demo Proof.

One first shows that if $\Omega_1$ is an open subset of $\IR^n$,
$\Omega_2$ is an open subset of $\IR^n$ 
with compact closure  contained in $\Omega_1$,
and $f\in L^p(\Omega_1)$ only depends on $(x_1,\ldots,x_d)$,
then $f\in L^p_d(\Omega_2)$ and
$\|f\|_{L^p_d(\Omega_2)} \le c \|f\|_{L^p(\Omega_1)}$.
The idea is that any ball of radius $r$ in $\Omega_2$
has the order of magnitude $r^{d-n}$ disjoint translates in $\Omega_1$
in the $(x_{d+1},\ldots,x_n)$ directions
and $f$ has the same $L^p$-norm on all these balls;
we leave the detailed verification to the reader.

Next consider the case of a compact manifold $X$ and $k=0$.
In a neighborhood of each point on $X$ there exist local coordinates
such that points with the same $(x_1,\ldots,x_d)$ lie in the same $H$-orbit.
The lemma then follows by applying the above estimate to these
coordinate charts.
For $k>0$ the lemma follows by applying the case $k=0$
to the functions $|\nablaA^js|$.
For $k<0$ the lemma follows by applying the case $k=0$
to the functions $|t_j|$,
once we observe that by averaging we may take the sections $t_j$
in (A.6) to be $H$-invariant.
\enddemo

\proclaim Lemma~A.2.

Multiplication gives bounded linear maps
$$L^p_\mi(X,E)\times L^q_\mi(X,F) \to L^r_\mi(X,E\otimes F)$$
for all $1\le p,q,r<\infty$ such that
$1/p+1/q=1/r$.

Using the norms (A.7)
we see that it suffices to show that multiplication gives bounded linear maps
$L^p_\mi(\IR^n)\times L^q_\mi(\IR^n)\to L^r_\mi(\IR^n)$.
That follows immediately from the H\"older inequality applied to $\IR^n$
and to the balls $B_\rho(x)$ in (A.1).

\proclaim Lemma~A.3.

Let $\mi\in(0,n)$.
Then there are continuous embeddings
$$L^{1,p}_\mi(X,E)\to  \cases{
              L^{p^*}_\mi(X,E) & for $p\in(1,\mi)$ \cr
               C^{0,\alpha}(X,E) & for $p\in(\mi,\infty)$ \cr
             }
$$
where $1/{p^*}=1/p-1/\mi$ and $\alpha=1-\mi/p$.

Again it suffices to establish the analogous
embedding on $\IR^n$.
The first  embedding 
is due to D.R.~Adams [1] Theorem~3.2; see also [4]~Theorem~2.
The second embedding 
is a classical result by C.B.~Morrey;
see [8] Chapter~3, Proposition~1.2 and Theorem~1.2.

\proclaim Lemma~A.4.

Let $\mi\in(0,n)$ and $p\in(1,\mi)$.
Let $q\in[1,p^*)$ where $1/p^*=1/p-1/\mi$.
Then the embedding $L^{1,p}_\mi(X,E)\to  L^q_\mi(X,E)$
is compact.

We have not found a reference for this particular result,
so we include a sketch of the proof.
It is of course similar to the proofs of other compact
embedding theorems.

\demo Proof.

It suffices to show that if $\Omega$
is a bounded open subset of $\IR^n$,
then the embedding $L^{1,p}_\mi(\Omega)\to  L^q_\mi(\Omega)$
is compact.
We may assume that $q\in[p,p^*)$.
Let $\alpha=\mi/q-\mi/p^*\in(0,1]$.
We say that
$f\in L^q_{\mi,\alpha}(\IR^n)$
if $f\in L^q_\mi(\IR^n)$ and
$$\sup_{0<|h|\le 1} |h|^{-\alpha}
         \bigl\| f(\cdot+h)-f(\cdot)\bigr\|_{L^q_\mi(\IR^n)} < \infty.
$$
This is a Banach space with norm
$$\|f\|_{L^q_{\mi,\alpha}(\IR^n)}^q
  = \|f\|_{L^q_\mi(\IR^n)}^q
    + \sup_{0<|h|\le 1} |h|^{-\alpha q}
            \bigl\| f(\cdot+h)-f(\cdot)\bigr\|_{L^q_\mi(\IR^n)}^q .
$$
Using the identity 
$f(x)=c_n^{-1}\sum_{j=1}^n\int_{\IR^n}
      (x_j-y_j)\,|x-y|^{-n}\,\partial_j f(y)\,dy
$
we get
$$\eqalign{
  & |h|^{-\alpha}\bigl(f(x+h)-f(x)\bigr) \cr
  & \quad = c_n^{-1} \sum_{j=1}^n \int_{\IR^n} |h|^{-\alpha} \bigl(
      (x_j-y_j+h_j)\,|x-y+h|^{-n} - (x_j-y_j)\,|x-y|^{-n}
                                            \bigr) \partial_j f(y)\,dy . \cr
          }
$$
Now 
$|h|^{-\alpha}\bigl|(z_j+h_j)\,|z+h|^{-n} - z_j|z|^{-n}\bigr|
  \le c \bigl(|z+h|^{1-n-\alpha}+|z|^{1-n-\alpha}\bigr)
$,
so
$$|h|^{-\alpha}\bigl|f(x+h)-f(x)\bigr|
  \le c\int_{\IR^n}
      |x-y|^{1-n-\alpha} \bigl( |\nabla f(y)| + |\nabla f(y+h)| \bigr)\,dy .
$$
It then follows from [1] Theorem~3.1 or [4] Theorem~2
that the right hand side,
viewed as a function of $x$,
 is bounded in $L^q_\mi(\IR^n)$ uniformly in $h$.
We conclude that there is a continuous embedding
$L^{1,p}_\mi(\IR^n)\to L^q_{\mi,\alpha}(\IR^n)$.

We define $L^q_{\mi,\alpha}(\Omega)$
by extension to $\IR^n$ as in (A.3).
Then there is a continuous  embedding
$L^{1,p}_\mi(\Omega)\to L^q_{\mi,\alpha}(\Omega)$.
It is straightforward to show that if $\Omega$ is bounded,
then the embedding
$L^q_{\mi,\alpha}(\Omega)\to L^q_\mi(\Omega)$ is compact.
\enddemo

\proclaim Lemma~A.5.

Let $p\in(1,\infty)$ and $k\in\IZ$.
Then any elliptic partial differential operator
$L^{k+m,p}_\mi(X,E)\to L^{k,p}_\mi(X,F)$
of order $m$ with smooth coefficients is a Fredholm operator.

This follows from the estimates for singular
integral operators on Morrey spaces due to J.~Peetre [15] Theorem~1.1;
see also [4] Theorem~3 and [19] Proposition~3.3.

{\it Remark A.6.}
Elliptic boundary value problems have not been studied in Morrey spaces.
However, it is clear that the classical treatment of the Neumann and Dirichlet
problems for the Laplacian by even and odd reflection
across the boundary carries over to Morrey spaces.

{\it Remark A.7.}
In this paper we have worked with smooth connections
and smooth gauge transformations.
One can of course complete the space of connections
and the group of gauge transformations in $L^{1,2}_3$
and $L^{2,2}_3$ respectively.
It can be shown that
$C^\infty_0(\IR^n)$ is not dense in $L^{k,p}_\mi(\IR^n)$
for $\mi\in[0,n)$.
The closure of $C^\infty_0(\IR^n)$ in $L^{k,p}_\mi(\IR^n)$
for $\mi\in[0,n)$
is the space $L^{k,p}_{\mi,0}(\IR^n)$
of all $f\in L^{k,p}_\mi(\IR^n)$
such that $f(\cdot+h)\to f(\cdot)$ in $L^{k,p}_\mi(\IR^n)$
as $h\to0\in\IR^n$.
It is of course a Banach space with the same norm as $L^{k,p}_\mi(\IR^n)$.
By the usual techniques, one can define
$L^{k,p}_{\mi,0}(X,E)$,
which then is the closure of $C^\infty(X,E)$ in $L^{k,p}_\mi(X,E)$.
Thus the completions of $\eusA$ and $\eusG$ are the space
of connections in $L^{1,2}_{3,0}$
and the group of gauge transformations in $L^{2,2}_{3,0}$.

}


{\ninepoint

\heading \tenbf Appendix~B. Homomorphisms of compact Lie groups

\noindent
The purpose of this appendix is to prove the following lemma:

\proclaim Lemma~B.1.

If $H$ is a semisimple compact Lie group
and $K$ is a compact Lie group,
then there exist only finitely many $K$-conjugacy classes
of smooth homomorphisms $H\to K$.

This follows from the next two lemmas.
For $H$ and $K$ compact Lie groups,
we let $\Hom(H,K)$ denote the set of continuous,
and hence smooth,
homomorphisms $H\to K$.
We view $\Hom(H,K)$ as a subset of the Banach manifold $C(H,K)$
of continuous maps $H\to K$.
The group $K$ acts on $\Hom(H,K)$ and $C(H,K)$ by conjugation.

\proclaim Lemma~B.2.

If $H$ is a semisimple compact Lie group
and $K$ is a compact Lie group,
then $\Hom(H,K)$ is a compact subset of $C(H,K)$.

\demo Proof.

Let $\eufh$ and $\eufk$ be the Lie algebras of $H$ and $K$.
Let $\hom(\eufh,\eufk)$ be the space of linear maps
$\eufh\to\eufk$ and $\eufhom(\eufh,\eufk)$ the
set of Lie algebra homomorphisms $\eufh\to\eufk$.
Since $\eufh$ is semisimple,
any homomorphism $\eufh\to\eufk$ takes values in the semisimple
part $\eufk_\ss$ of $\eufk=\eufz(\eufk)\oplus\eufk_\ss$.
The Killing form $|X|^2=-\Tr(\ad X)^2$
defines norms on $\eufh$ and $\eufk_\ss$.
This gives a norm
$$|\lambda|^2  = \sup_{0\ne X\in\eufh} {|\lambda X|^2  \over |X|^2 }
               = \sup_{0\ne X\in\eufh}
                    { \Tr(\ad(\lambda X))^2 \over \Tr (\ad X)^2}
$$
on $\hom(\eufh,\eufk_\ss)$.
Now $\ad\circ\lambda$ gives a representation of $\eufh$ on $\eufk$.
We see that $|\lambda|$ only depends on
(the isomorphism class of) this representation.
It follows from the classification of 
representations of semisimple Lie algebras
that there are only finitely many representations of $\eufh$
of given dimension.
Hence $|\lambda|$ can assume only finitely
many values.
In particular, $\eufhom(\eufh,\eufk)$ is a bounded subset
of $\hom(\eufh,\eufk_\ss)$.

If $\sigma\in\Hom(H,K)$,
then $(\sigma_*)_1\in\eufhom(\eufh,\eufk)$.
Thus we have a uniform bound for $(\sigma_*)_1$ for all $\sigma\in \Hom(H,K)$.
By left invariance, 
this gives a uniform bound for $(\sigma_*)_h$ for all $\sigma\in\Hom(H,K)$
and $h\in H$.
Finally, $\Hom(H,K)$ is a closed subset of $C(H,K)$.
It then follows from the Arzela-Ascoli theorem that
$\Hom(H,K)$ is a compact subset of $C(H,K)$.
\enddemo

\proclaim Lemma~B.3. \rm(R.W.~Richardson [13], D.H.~Lee [12])

If $H$ and $K$are compact Lie groups,
then $\Hom(H,K)$ is a discrete union of $K$-orbits in $C(H,K)$.

We include a simplified version of the proof in [12].

\demo Proof.

We have to  show that any $K$-orbit in $\Hom(H,K)$ has a tubular neighborhood
in $C(H,K)$ that does not intersect any other $K$-orbits in $\Hom(H,K)$.
We write $1$ for the identity elements in $H$ and $K$.
We also write 1 for the map $H\times H\to K$
that maps $H\times H$ to 1.
Then $\Hom(H,K)=T^{-1}(1)$ where the map
$$\displaylines{
  T:C(H,K)\to C(H\times H,K) \cr
  \noalign{\line{is defined as\hfill}}
  T(\sigma)(h_1,h_2)=\sigma(h_1)\,\sigma(h_2)\,\sigma(h_1h_2)^{-1}. \cr
               }
$$
Formally the tangent space of $C(H,K)$ at $\sigma$ is given by
the null space of the differential $(T_*)_\sigma$.

Let $\sigma\in\Hom(H,K)$.
Then the $K$-orbit of $\sigma$ is given by the range
of the map
$$\displaylines{
  S:K\to C(H,K) \cr
  \noalign{\line{defined as\hfill}}
  S(k)(h)=k\,\sigma(h)\,k^{-1}. \cr
               }
$$
This orbit is a homogeneous $K$-space.
Its tangent space at $\sigma$ to the $K$-orbit of $\sigma$
is given by the range of the differential $(S_*)_1$.

We identify the tangent space at any point in the Lie group $K$
with the Lie algebra $\eufk$ by right translation.
For any compact topological space $\eufX$
this gives an identification of the tangent space
at any point in the Banach manifold $C(\eufX,K)$ 
with the Banach space $C(\eufX,\eufk)$.
Under these identifications, the differentials
$$\eqalignno{
    (S_*)_1&:\eufk\to C(H,\eufk) \cr
    (T_*)_\sigma&:C(H,\eufk)\to C(H\times H,\eufk) \cr
\noalign{\line{are given by\hfil}}
    (S_*)_1(\xi)(h)&=\xi-\Ad\sigma(h)\,\xi \cr
    (T_*)_\sigma(\lambda)(h_1,h_2)&=\lambda(h_1)+\Ad\sigma(h_1)\,\lambda(h_2)
        -\lambda(h_1h_2) . \cr
               }
$$
We see that $(S_*)_1=-\delta_1$ and $(T_*)_\sigma=\delta_2$ where
the linear maps
$$\delta_q:C(H^{q-1},\eufk) \to C(H^q,\eufk)$$
are defined as
$$\eqalign{
  \delta_q\mu(h_1,\ldots,h_q) & = \Ad\sigma(h_1)\,\mu(h_2,\ldots,h_q) \cr
    & \qquad  + \sum_{i=1}^{q-1} (-1)^i \mu(h_1,\ldots,h_ih_{i+1},\ldots,h_q)
              + (-1)^q \mu(h_1,\ldots,h_{q-1}) . \cr
          }
$$
A short calculation shows that the linear maps $\delta_q$ define
a chain complex.
Its cohomology groups $H^q_\cont(H,\eufk)$
are known as the continuous cohomology groups
of $H$ with coefficients in the $H$-module $\eufk$.
We have assumed $H$ to be compact,
so we can define maps
$$\displaylines{
    s_q: C(H^q,\eufk) \to C(H^{q-1},\eufk) \cr
   \noalign{\line{by\hfill}}
    s_q \mu (h_1,\ldots,h_{q-1})
       = (-1)^q \int_H \mu (h_1,\ldots,h_{q-1},h)\,dh , \cr
               }
$$
where $dh$ is the normalized Haar measure on $H$.
Then a short calculation shows that
$$\cases{
  s_1\delta_1 = 1-p \cr
  \delta_q s_q+s_{q+1}\delta_{q+1} = 1 & for $q\ge1$, \cr
        }
$$
where $p(\xi)=\int_H \Ad\sigma(h)\,\xi\,dh$.
The map $p$ is a projection of
$\eufk$ onto the fixed point set $\eufk^H$.
It follows that
$$H^q_\cont(H,\eufk) = \cases{ \eufk^H & for $q=0$ \cr 0 & for $q\ge1$ .}$$
That $H^1_\cont(H,\eufk)=0$ tells us that the null space of $(T_*)_\sigma$
is the tangent space of the $K$-orbit of $\sigma$.
That $H^2_\cont(H,\eufk)=0$ tells us that $(T_*)_\sigma$ has closed range.
It then follows from the implicit function theorem that
the $K$-orbit of $\sigma$ has a tubular neighborhood in $C(H,K)$
where $T\ne1$ away from the $K$-orbit itself.
\enddemo

}

Continuous cohomology of Lie groups was introduced in [5].
For modern treatments, see [11] and the references listed there.
Continuous cohomology is essentially group cohomology
with continuous cochains.
It is the natural cohomology theory for Lie groups.


\headingfeedfalse
\references

\ref[1]
D.R. Adams,
{\it A note on Riesz potentials},
Duke Math. J. {\bf42} (1975), 765--778.

\ref[2]
G. Bor,
{\it Yang-Mills fields which are not self-dual},
Comm. Math. Phys. {\bf145} (1992), 393--410.

\ref[3]
P.J. Braam and G. Mati\'c,
{\it The Smith conjecture in dimension four and equivariant gauge theory},
Forum Math. {\bf5} (1993), 299--311.

\ref[4]
F. Chiarenza and M. Frasca,
{\it Morrey spaces and Hardy-Littlewood maximal function},
Rend. Mat. Appl. (7) {\bf7} (1987), 273--279.

\ref[5]
W.T. van Est,
{\it On the algebraic cohomology concepts in Lie groups I and II},
Indag. Math. {\bf17} (1955), 225--233, 286--294.

\ref[6]
R. Fintushel and R. Stern,
{\it Pseudofree orbifolds},
Ann. of Math. {\bf122} (1985), 335--364.

\ref[7]
M. Furuta,
{\it A remark on a fixed point of finite group action on $S^4$},
Topology {\bf28} (1989), 35--38.

\ref[8]
M. Giaquinta,
{\it Multiple integrals in the calculus of variations and nonlinear
elliptic systems},
Princeton Univ. Press, Princeton, 1983.

\ref[9]
U. Gritsch,
{\it Morse theory for the Yang-Mills functional over equivariant
four-manifolds and equivariant homotopy theory},
Ph.D. thesis, Stanford University, 1997.

\ref[10]
\vrule height 0pt depth 0.3pt width 20pt \thinspace ,
{\it Morse theory for the Yang-Mills functional 
via equivariant homotopy theory},
to appear in Trans. Amer. Math. Soc.

\ref[11]
A. Guichardet,
{\it Cohomologie des groupes topologiques et des alg\`ebres de Lie},
Cedic-Nathan, Paris, 1980.

\ref[12]
D.H. Lee,
{\it On deformations of homomorphisms of locally compact groups},
Trans. Amer. Math. Soc. {\bf 191} (1974), 353--361.

\ref[13]
A. Nijenhuis and R.W. Richardson, Jr.,
{\it Deformations of homomorphisms of Lie groups and Lie algebras},
Bull. Amer. Math. Soc. {\bf73} (1967), 175--179.

\ref[14]
T.H. Parker,
{\it A Morse theory for equivariant Yang-Mills},
Duke Math. J. {\bf66} (1992), 337--356.

\ref [15]
J. Peetre,
{\it On convolution operators leaving $L^{p,\lambda}$ spaces invariant},
Ann. Mat. Pura Appl. (4) {\bf 72} (1966), 295--304.

\ref[16]
L. Sadun and J. Segert,
{\it Non-self-dual Yang-Mills connections with quadrupole symmetry},
Comm. Math. Phys. {\bf145} (1992), 363--391.

\ref[17]
S. Sedlacek,
{\it A direct method for minimizing the Yang-Mills functional},
Commun. Math. Phys. {\bf86} (1982), 515--527.

\ref[18]
C.H. Taubes,
{\it Path-connected Yang-Mills moduli spaces},
J. Differential Geom. {\bf19} (1984), 337--392.

\ref [19]
M.E. Taylor,
{\it Analysis on Morrey spaces and applications to Navier-Stokes and
other evolution equations},
Comm. Partial Differential Equations {\bf17} (1992), 1407--1456.

\ref[20]
K.K.~Uhlenbeck,
{\it Connections with $L^p$ bounds on curvature},
Comm. Math. Phys. {\bf83} (1982), 31--42.

\ref [21]
H. Urakawa,
{\it Equivariant theory of Yang-Mills connections over Riemannian
manifolds of cohomogeneity one},
Indiana Univ. Math. J. {\bf37} (1988), 753--788.

\ref[22]
H.-Y. Wang,
{\it The existence of nonminimal solutions to the Yang-Mills equation
with group $\SU(2)$ on $S^2\times S^2$ and $S^1\times S^3$},
J. Differential Geom. {\bf34} (1991), 701--767.

\ref[23]
D. Wilkins,
{\it The Palais-Smale conditions for the Yang-Mills functional},
Proc. Roy. Soc. Edinburgh Sect. A {\bf108} (1988), 189-200.

\endlines{Department of Mathematics, G\"oteborg University,
  S-412\thinspace96 G\"oteborg, Sweden\cr
  http://www.math.chalmers.se/$\sim$jrade/\cr
  jrade@math.chalmers.se\cr 
  April 7, 2000\cr
         }
 
\bye